\documentclass[11pt]{amsart}
\usepackage{amscd,amssymb,amsmath}
\usepackage[arrow,matrix]{xy}
\usepackage{graphicx}
\usepackage{cite}
\newtheorem{thm}[subsection]{Theorem}

\newtheorem{proposition}[subsection]{Proposition}
\newtheorem{cor}[subsection]{Corollary}

\newtheorem{rk}[subsection]{Remark}
\newtheorem{defn}[subsection]{Definition}

\numberwithin{equation}{section} 

\newcommand{\B}{{\mathcal B}}

\newcommand{\bea}{\begin{eqnarray}}
\newcommand{\eea}{\end{eqnarray}}

\newcommand{\R}{\mathbb{R}}

\DeclareMathOperator{\Der}{\rm Der}

\begin{document}
\title [Evolution algebra of a bisexual population]
{Evolution algebra of a bisexual population}

\author {M. Ladra, U.A. Rozikov}

 \address{M. Ladra\\ Department of Algebra, University of Santiago de Compostela, 15782, Spain.}
 \email {manuel.ladra@usc.es}
 \address{U.\ A.\ Rozikov\\ Institute of mathematics and information technologies,
Tashkent, Uzbekistan.}
\email {rozikovu@yandex.ru}

\begin{abstract} We introduce an (evolution) algebra identifying the coefficients
of inheritance of a bisexual population as the structure constants
of the algebra. The basic properties of the algebra are studied. We
prove that this algebra is commutative (and hence flexible), not
associative and not necessarily power associative. We show that the
evolution algebra of the bisexual population is not a baric algebra,
but a dibaric algebra and hence its square is baric. Moreover, we
show that the algebra is a Banach algebra. The set of all
derivations of the evolution algebra is described. We find necessary
conditions for a state of the population to be a fixed point or a
zero point of the evolution operator which corresponds to the
evolution algebra. We also establish upper estimate of the limit
points set for trajectories of the evolution operator. Using the
necessary conditions we give a detailed analysis of a special case
of the evolution algebra (bisexual population of which has a
preference on type ``1'' of females and males). For such a special
case we describe the full set of idempotent elements and the full
set of absolute nilpotent elements.

\end{abstract}
\maketitle
\section{Introduction} \label{sec:intro}

The action of genes is manifested statistically in sufficiently
large communities of matching individuals (belonging to the same
species). These communities are called {\it populations} \cite{ly}.
The population exists not only in space but also in time, i.e. it
has its own life cycle. The basis for this phenomenon is
reproduction by mating. Mating in a population can be free or
subject to certain restrictions.

The whole population in space and time comprises discrete
generations $F_0, F_1, \dots$ The generation $F_{n+1}$ is the set of
individuals whose parents belong to the $F_n$ generation. A state of
a population is a distribution of probabilities of the different
types of organisms in every generation. Type partition is called
differentiation. The simplest example is sex differentiation. In
bisexual population any kind of differentiation must agree with the
sex differentiation, i.e. all the organisms of one type must belong
to the same sex. Thus, it is possible to speak of male and female
types.

The evolution (or dynamics) of a population comprises a determined
change of state in the next generations as a result of reproduction
and selection. This evolution of a population can be studied by a
dynamical system (iterations) of a quadratic stochastic operator.

The history of the quadratic stochastic operators can be traced back
to the work of S. Bernstein\cite{ber}. During more than 85 years
this theory developed and many papers were published (see e.g.
\cite{ber,gr,gan1,ge,gar,ly,rz,rz1,roz,rs}). In recent years it has
again become of interest in connection with numerous applications to
many branches of mathematics, biology and physics.

A quadratic stochastic operator (QSO), $V$, has meaning of a {\it
free} population evolution operator, which arises as follows:
Consider a free population consisting of $m$ species. Let
$x^0=(x_1^0,\dots,x_m^0)$ be the probability distribution of species
in the initial generations, and $P_{ij,k}$ the probability that
individuals in the $i$th and $j$th species interbreed to produce an
individual $k$. Then the probability distribution
$x'=(x_1',\dots,x'_m)$ of the species in the first generation can be
found by the total probability i.e.

\begin{equation}\label{1}x'_k=\sum_{i,j=1}^{m}P_{ij,k}x_i^0x_j^0,\ k=1,\dots,m\end{equation}
 where the cubic matrix $P\equiv P(V)=(P_{ij,k})_{i,j,k=1}^m$ satisfies the following conditions

 \begin{equation}\label{2}P_{ij,k}\geq 0,\ \ \sum_{k=1}^m P_{ij,k}=1,\ i,j\in \{1,\dots,m\}.
\end{equation}

 This means that the association $x^0\mapsto x'$ defines a map $V$ of the simplex
  \begin{equation}S^{m-1}=\{x=(x_1,\dots,x_m)\in \R^m\ :\ x_i\geq 0, \sum_{i=1}^m x_i=1\} \end{equation}
  into itself, called the {\it evolution operator}.

The population evolves by starting from an arbitrary state
$x^{(0)}$, then passing to the state
$x'=V(x^{(0)})$, then to the state $x''=V(V(x^{(0)})), $ and so on.\\

There are many papers devoted to study of the evolution of the free
population, i.e. to study of dynamical system generated by the QSO
(\ref{1}), see e.g. \cite{gr,gan1,ge,gar,ly,rz,rz1,roz,rs}.

There exist several classes of non-associative algebras (baric,
evolution, Bernstein, train, stochastic, etc.), whose investigation
has provided a number of significant contributions to theoretical
population genetics. Such classes have been defined different times
by several authors, and all algebras belonging to these classes are
generally called "genetic." Etherington introduced the formal
language of abstract algebra to the study of the genetics in his
series of seminal papers \cite{e1,e2,e3}.  In recent years many
authors have tried to investigate the difficult problem of
classification of these algebras. The most comprehensive references
for the mathematical research done in this area are \cite{ly,m,t,w}.

 In \cite{ly} an evolution algebra $\mathcal A$ associated to the
free population is introduced and using this non-associative algebra
many results are obtained in explicit form, e.g. the explicit
description of stationary quadratic operators, and the explicit
solutions of a nonlinear evolutionary equation in the absence of
selection, as well as general theorems on convergence to equilibrium
in the presence of selection.

In \cite{t} a new type of evolution algebra is introduced. This
algebra also describes some evolution laws of genetics and it is an
algebra $E$ over a field $K$ with a countable natural basis
$e_1,e_2,\dots$ and multiplication given by $e_ie_i=\sum_j
a_{ij}e_j$, $e_ie_j=0$ if $i\ne j$. Therefore, $e_ie_i$ is viewed as
``self-reproduction''.

 In this paper we consider a bisexual population (BP) and define an evolution
 algebra (EA) using inheritance coefficients of the population. This algebra is a natural
 generalization of the algebra $\mathcal A$ of free population. The evolution algebra of a bisexual
 population (EABP) is different from the EA defined in \cite{t}.

The paper is organized as follows. In Section 2 we give evolution
operator of BP. Section 3 contains the definition of EABP. Section 4
is devoted to basic properties of the algebra and therein we prove
that the EABP is commutative but not associative and not power
associative. In Section 5 we show that the EABP is not a baric
algebra. In Section 6 we prove that the EABP is a dibaric algebra,
hence its square is a baric algebra. In Section 7 we prove that the
EABP is a Banach algebra. In Section 8 we describe the set of all
derivations of EABP. Section 9 is devoted to study dynamics of the
evolution operator, which corresponds to the evolution algebra. We
find necessary conditions for a state of the population to be a
fixed point or a zero point of the evolution operator. We also
establish upper estimate of the limit points set for trajectories of
the evolution operator. In the last section we give a detailed
analysis of a special case of the evolution algebra (bisexual
population of which has a preference on type ``1'' of females and
males). For such a special case we describe the full set of
idempotent elements and the full set of absolute nilpotent elements.

\section{Evolution operator of a BP}

  In this section following \cite{ly}, we describe the evolution operator of a BP.
Assuming that the population is bisexual we suppose that the set of
females can be partitioned into finitely many different types
indexed by $\{1,2,\dots,n\}$ and, similarly, that the male types are
indexed by $\{1,2,\dots,\nu \}$. The number $n+\nu$ is called the
dimension of the population. The population is described by its
state vector $(x,y)$ in $S^{n-1}\times S^{\nu-1}$, the product of
two unit simplexes   in $\R^n$ and $\R^\nu$ respectively. Vectors
$x$ and $y$ are the probability distributions of the females and
males over the possible types:
$$x_i\geq 0, \ \sum_{i=1}^{n}x_i=1; \ \ y_i\geq 0, \
\sum_{i=1}^{\nu}y_i=1.$$

Denote $S=S^{n-1}\times S^{\nu-1}$. We call the partition into types
hereditary if for each possible state $z=(x,y)\in S$ describing the
current generation, the state $z'=(x',y')\in S$ is uniquely defined
describing the next generation. This means that the association
$z\mapsto z'$ defines a map $V:S\to S$ called the evolution
operator.

For any point $z^{(0)}\in S$ the sequence $z^{(t)}=V(z^{(t-1)}),
t=1,2,\dots$ is called the trajectory of $z^{(0)}$.

 Let $P_{ik,j}^{(f)}$ and $P_{ik,l}^{(m)}$ be inheritance coefficients
 defined as the probability that a female offspring is type $j$ and, respectively,
 that a male offspring is of type $l$, when the parental pair is
 $ik$ $(i,j=1,\dots,n$; and $k,l=1,\dots,\nu)$. We have
\begin{equation}\label{2}
P_{ik,j}^{(f)}\geq 0, \ \ \sum_{j=1}^nP_{ik,j}^{(f)}=1; \ \
P_{ik,l}^{(m)}\geq 0, \ \ \sum_{l=1}^\nu P_{ik,l}^{(m)}=1.
\end{equation}

Let $z'=(x',y')$ be the state of the offspring population at the
birth stage. This is obtained from inheritance coefficients as

\begin{equation}\label{3}
x'_j= \sum_{i,k=1}^{n,\nu}P_{ik,j}^{(f)}x_iy_k; \ \ y'_l=
\sum_{i,k=1}^{n,\nu} P_{ik,l}^{(m)}x_iy_k.
\end{equation}

We see from (\ref{3}) that for a BP the evolution operator is a
quadratic mapping of $S$ into itself. But for free population the
operator is quadratic mapping of the simplex into itself given by
(\ref{1}).

\section{Definition of the EABP}
 In this section we give an algebra structure on the vector space
 $\R^{n+\nu}$ which is closely related to the map (\ref{3}).

 Consider $\{e_1,\dots,e_{n+\nu}\}$ the canonical basis on $\R^{n+\nu}$
 and divide the basis as $e^{(f)}_i=e_i$, $ i=1,\dots,n$ and $e^{(m)}_i=e_{n+i}$,
 $i=1,\dots,\nu$.

 Now introduce on $\R^{n+\nu}$ a multiplication defined by
\begin{equation}\label{4}
\begin{array}{ll}
e^{(f)}_ie^{(m)}_k = e^{(m)}_ke^{(f)}_i={1\over 2}\left(\sum_{j=1}^nP_{ik,j}^{(f)}e^{(f)}_j+ \sum_{l=1}^{\nu}P_{ik,l}^{(m)}e^{(m)}_l\right),  \\
e^{(f)}_ie^{(f)}_j =0,\ \ i,j=1,\dots,n; \ \ e^{(m)}_ke^{(m)}_l =0,
\ \ k,l=1,\dots,\nu. \end{array} \end{equation}

Thus we identify the coefficients of bisexual inheritance as the
structure constants of an algebra, i.e. a bilinear mapping of
$\R^{n+\nu}\times \R^{n+\nu}$ to $\R^{n+\nu}$.

The general formula for the multiplication is the extension of
(\ref{4}) by bilinearity, i.e. for $z,t\in \in \R^{n+\nu}$,
$$z=(x,y)=\sum_{i=1}^nx_ie_i^{(f)}+\sum_{j=1}^\nu y_je_j^{(m)}, \ \ t=(u,v)=\sum_{i=1}^nu_ie_i^{(f)}+\sum_{j=1}^\nu v_je_j^{(m)}$$
using (\ref{4}), we obtain
\begin{equation}\label{5}
\begin{array}{ll} zt={1\over 2}\sum_{k=1}^n\left(\sum_{i=1}^n\sum_{j=1}^\nu
P_{ij,k}^{(f)}(x_iv_j+u_iy_j)\right)e^{(f)}_k+ \\[2mm] {1\over
2}\sum_{l=1}^\nu\left(\sum_{i=1}^n\sum_{j=1}^\nu
P_{ij,l}^{(m)}(x_iv_j+u_iy_j)\right)e^{(m)}_l.\end{array}
\end{equation}

From (\ref{5}) and using (\ref{3}), in the particular case that
$z=t$, i.e. $x=u$ and $y=v$, we obtain
\begin{equation}\label{6}
\begin{array}{ll} zz=z^2=\sum_{k=1}^n\left(\sum_{i=1}^n\sum_{j=1}^\nu
P_{ij,k}^{(f)}x_iy_j\right)e^{(f)}_k+ \\[2mm]
\sum_{l=1}^\nu\left(\sum_{i=1}^n\sum_{j=1}^\nu
P_{ij,l}^{(m)}x_iy_j\right)e^{(m)}_l=V(z)\end{array}
\end{equation}
for any $z\in S$.

 This algebraic interpretation is very useful. For
example, a BP state $z=(x,y)$ is an equilibrium (fixed point,
$V(z)=z$) precisely when $z$ is an idempotent element of the set
$S$.

If we write $z^{[t]}$ for the power $(\cdots(z^2)^2\cdots)$ ($t$
times) with $z^{[0]}\equiv z$ then the trajectory with initial state
$z$ is $V^t(z)=z^{[t]}$.

The algebra ${\mathcal B}={\mathcal B}_V$ generated by the evolution
operator $V$ (see (\ref{3})) is called the {\it evolution algebra of
the bisexual population} (EABP).

  \begin{rk}
{\rm 1.} If a population is free then the male and female types are
identical and, in particular $n=\nu$, the inheritance coefficients
are the same for male and female offsprings, i.e.
$$P_{ik,j}= P_{ik,j}^{(f)}=P_{ik,j}^{(m)}.$$

The evolution algebra $\mathcal A$ associated with the free
population is introduced and studied in \cite{ly}. Note that this
algebra is commutative when the condition of symmetry
$P_{ik,j}=P_{ki,j}$ is satisfied, but it is not in general
associative. In the next section we show that algebra $\mathcal B$
of bisexual population is commutative without any symmetry
condition. Hence the algebra $\mathcal A$ is a particular case of
the algebra $\mathcal B$.

{\rm 2}. It is easy to see that the EA introduced in \cite{t} is
completely different from our EABP, i.e. of $\mathcal B$.

{\rm 3}. The algebra $\mathcal B$ is a natural generalization of a
{\it zygotic algebra} for sex linked inheritance (see
\cite{e3,h,m,w}).
\end{rk}
\section{ Basic properties of the EABP}

The following theorem gives basic properties of the EABP.

\begin{thm}\label{t1} 1) Algebra $\mathcal B$ is not associative, in
general.

2) Algebra $\mathcal B$ is commutative, flexible.

3) $\mathcal B$ is not power-associative, in general.
\end{thm}
\proof 1) Take $e_i^{(f)}$, $e_j^{(m)}$ such that $P_{ij,s}^{(f)}\ne
0$ for some $s$ and take $e_k^{(m)}$ such that $P_{sk,r}^{(f)}\ne 0$
for some $r$ then
$$(e_i^{(f)}e_j^{(m)})e_k^{(m)}={1\over 2}\sum_{q=1}^nP_{ij,q}^{(f)}e_q^{(f)}e_k^{(m)}=$$
$${1\over 2}\left(P_{ij,s}^{(f)}e_s^{(f)}e_k^{(m)}+ \sum_{q=1\atop q\ne s}^nP_{ij,q}^{(f)}e_q^{(f)}e_k^{(m)}\right)=$$
$${1\over 4}P_{ij,s}^{(f)}P_{sk,r}^{(f)}e_r^{(f)}+\mbox{non-negative terms} \ne 0.$$
But $$ e_i^{(f)}(e_j^{(m)}e_k^{(m)})=0, \ \ \mbox{i.e.} \
\ (e_i^{(f)}e_j^{(m)})e_k^{(m)}\ne e_i^{(f)}(e_j^{(m)}e_k^{(m)}). $$

2) Commutativity of $\mathcal B$ follows from formula (\ref{5}). An
algebra is called {\it flexible} if it satisfies $z(tz)=(zt)z$ for
any $z,t$. It is easy to see that a commutative algebra is flexible.

3) To show that $\mathcal B$ is not a power-associative, in general,
we shall construct an example of $z$ such that $(zz)(zz)\ne
((zz)z)z.$ Consider $n=1$, $\nu=2$. In this case
$$P_{11,1}^{(f)}=P_{12,1}^{(f)}=1, \ \
P_{11,1}^{(m)}+P_{11,2}^{(m)}=1, \ \
P_{12,1}^{(m)}+P_{12,2}^{(m)}=1.$$ Denote $a=P_{11,1}^{(m)}, \ \
b=P_{12,1}^{(m)}$. Take $z=e_1^{(f)}+e_1^{(m)}$. Then we have
$$z^2=e_1^{(f)}+a e_1^{(m)}+(1-a)e_2^{(m)}.$$
\begin{equation}\label{7}
z^2z^2=e_1^{(f)}+\left(a^2+(1-a)b\right)e_1^{(m)}+(1-a)(a-b+1)e_2^{(m)}.
\end{equation}
$$z^2z=e_1^{(f)}+{1\over 2}\left(a^2+(1-a)b+a\right)e_1^{(m)}+{1\over
2}(1-a)(a-b+2)e_2^{(m)}.$$
\begin{equation}\label{8}\begin{array}{ll}
(z^2z)z=e_1^{(f)}+{1\over
4}\left(a(a-b)^2+(a+b)(a-b+2)\right)e_1^{(m)}+\\[2mm]
\ \ \ \ \ \ \ \ \ \ \ {1\over
4}(1-a)\left(3+(a-b+1)^2\right)e_2^{(m)}.\end{array}
\end{equation}
Assume $a=P_{11,1}^{(m)}\ne 1$ and $a\ne b=P_{12,1}^{(m)}$. Then
from (\ref{7}) and (\ref{8}) we get $(zz)(zz)\ne ((zz)z)z.$
\endproof

\section{ $\mathcal B$ is not a baric algebra}

A {\it character} for an algebra $A$ is a nonzero multiplicative
linear form on $A$, that is, a nonzero algebra homomorphism from $A$
to $\R$ \cite{ly}. Not every algebra admits a character. For
example, an algebra with the zero multiplication has no character. A
pair $(A, \sigma)$ consisting of an algebra $A$ and a character
$\sigma$ on $A$ is called a {\it baric algebra}. In \cite{ly} for
the EA of a free population it is proven that there is a character
$\sigma(x)=\sum_i x_i$, therefore that algebra is baric. But the
following theorem says that the EABP, i.e. $\mathcal B$ is not
baric.

\begin{thm}\label{t2} The EABP, $\mathcal B$, has no a nonzero character.
\end{thm}
\proof Assume $\sigma(z)=\sum_{i=1}^n a_i x_i+ \sum_{j=1}^\nu b_j
y_j$, $ z=(x,y)\in \mathcal B$ is a character. We shall prove that
$\sigma(z)\equiv 0$. For $z=(x,y)$, $t=(u,v)$ we have
$$\sigma(zt)={1\over
2}\sum_{i=1}^n\sum_{j=1}^\nu\left(\sum_{p=1}^n
a_pP_{ij,p}^{(f)}+\sum_{q=1}^\nu b_q
P_{ij,q}^{(m)}\right)(x_iv_j+u_iy_j);$$
$$\sigma(z)\sigma(t)=\sum_{i=1}^n\sum_{j=1}^n a_ia_jx_iu_j+\sum_{i=1}^n\sum_{j=1}^\nu
a_ib_j(x_iv_j+u_iy_j)+\sum_{i=1}^\nu\sum_{j=1}^\nu b_ib_jy_jv_j.$$
From $\sigma(zt)=\sigma(z)\sigma(t)$ we get
$$ a_ia_j=0\ \ \mbox{for any } \ \ i,j=1,\dots,n$$
$$ b_ib_j=0\ \ \mbox{for any } \ \ i,j=1,\dots,\nu$$
This system of equations has only the solution
$a_1=\dots=a_n=b_1=\dots=b_{\nu}=0$.
\endproof
 \begin{rk}\label{rr}
In \cite{ly} the baric EA of a free population is classified as
algebra induced by a linear operator; unit algebra; constant
algebra; Bernstein (stationary) algebra; genetic algebra; train
algebra, etc. But since the EABP is not baric there are not similar
algebras for bisexual population.
\end{rk}

\section{ $\mathcal B$ is a dibaric algebra}

 By Theorem \ref{t2} the algebra  $\mathcal B$ is not a baric
 algebra. To overcome such complication, Etherington \cite{e3} for
 a zygotic algebra of sex linked inheritance introduced the idea of
 treating the male and female components of a population separately.
 In \cite{h} Holgate formalized this concept by introducing sex
 differentiation algebras and a generalization of baric algebras called dibaric algebras.
 In this section we shall prove that the algebra  $\mathcal B$ is a
 dibaric algebra.

\begin{defn}\cite{m,w} Let $\mathfrak A=\langle w, m \rangle_\R$ denote
a two dimensional commutative algebra over $\R$ with multiplicative
table
$$w^2=m^2=0, \ \ wm={1\over 2}(w+m).$$
Then  $\mathfrak A$ is called the {\it sex differentiation algebra}.
\end{defn}

As usual, a {\it subalgebra} $B$ of an algebra $A$ is a subspace
which is closed under multiplication. A subspace $B$ is an {\it
ideal} if it is closed under multiplication by all elements in $B$.
For example, the {\it square} of the algebra:
$$A^2={\rm span}\{zt: z,t\in A\}$$ is an ideal.

 It is clear that
$\mathfrak A^2=\langle w+m \rangle_\R$ is an ideal of $\mathfrak A$
which is isomorphic to the field $\R$. Hence the algebra $\mathfrak
A^2$ is a baric algebra. Now we can define Holgate's generalization
of a baric algebra.

\begin{defn}\cite{m} An algebra is called dibaric if it admits a
homomorphism onto the sex differentiation algebra $\mathfrak A$.
\end{defn}

\begin{thm}\label{td} The algebra $\mathcal B$ is dibaric.
\end{thm}
\proof Consider mapping $\varphi:\B\to\mathfrak A$ defined by
\begin{equation}\label{mm} \varphi(e_i^{(f)})=w, \ \ i=1,\dots,n; \ \
\varphi(e_k^{(m)})=m, \ \ i=1,\dots,\nu.\end{equation}
 For $z,t\in \R^{n+\nu}$,
$$z=(x,y)=\sum_{i=1}^nx_ie_i^{(f)}+\sum_{j=1}^\nu y_je_j^{(m)}, \ \
t=(u,v)=\sum_{i=1}^nu_ie_i^{(f)}+\sum_{j=1}^\nu v_je_j^{(m)},$$
using linearity of $\varphi$, (\ref{mm}), (\ref{2}) and (\ref{5}) we
get
$$\varphi(zt)=\sum_{i=1}^n\sum_{j=1}^\nu\left(x_iv_j+u_iy_j\right)\left({1\over
2}(w+m)\right).$$ We also have
$$\varphi(z)\varphi(t)=\left(\sum_{i=1}^nx_iw+\sum_{j=1}^\nu y_jm\right)
\left(\sum_{i=1}^nu_iw+\sum_{j=1}^\nu v_jm\right)=$$
$$\sum_{i=1}^n\sum_{j=1}^\nu\left(x_iv_j+u_iy_j\right)\left({1\over
2}(w+m)\right)=\varphi(zt),$$ i.e. $\varphi$ is a homomorphism. For
arbitrary $u=\alpha w+\beta m\in \mathfrak A$ it is easy to see that
$\varphi(z)=u$ if $z=\sum_{i=1}^nx_ie_i^{(f)}+\sum_{j=1}^\nu
y_je_j^{(m)}\in \B$ with $\sum_{i=1}^nx_i=\alpha$ and
$\sum_{j=1}^\nu y_j=\beta$. Therefore $\varphi$ is onto.
\endproof
\begin{proposition}\label{h}\cite{h} If an algebra $A$ is dibaric, then $A^2$
is baric.
\end{proposition}
As a corollary of this Proposition and Theorem \ref{td} we have
\begin{cor}\label{c} The subalgebra $\B^2$
is a baric algebra.
\end{cor}
Since $\B^2$ is a baric algebra all types of particular algebras
mentioned in Remark \ref{rr} can be defined and studied for the
algebra $\B^2$.

\section{$\B$ is a Banach algebra}

Define a norm $\|\cdot\|:\mathcal B\to \R$ as follows
$$\|z\|=\sum_{i=1}^n|x_i|+\sum_{j=1}^\nu|y_j|,$$
where $z=\sum_{i=1}^nx_ie^{(f)}_i+\sum_{j=1}^\nu y_j e^{(m)}_j\in
\mathcal B$.

For a fixed $a\in \mathcal B$ consider the operator $L_a: \mathcal
B\to \mathcal B$, left (right) multiplication, defined as
$$ L_a(z)=az \ \ (=za).$$
It is easy to see that the set $\left\{L_{e^{(f)}_i}, L_{e^{(m)}_j},
i=1,\dots,n, \ j=1,..,\nu\right\}$ spans a linear space which is the
set of all the operators of the left (right) multiplications.
\begin{proposition}
Operator $L_a$ is a bounded linear operator for any
$a=\sum_{i=1}^na_ie^{(f)}_i+\sum_{j=1}^\nu b_j e^{(m)}_j\in \mathcal
B$.
\end{proposition}
\proof For $z=\sum_{i=1}^nx_ie^{(f)}_i+\sum_{j=1}^\nu y_j
e^{(m)}_j\in \mathcal B$ we have
$$L_a(z)=\sum_{p=1}^n\left({1\over 2}\sum_{i=1}^n\sum_{j=1}^\nu
P_{ij,p}^{(f)}(a_iy_j+x_ib_j)\right)e^{(f)}_p+$$ $$
\sum_{q=1}^\nu\left({1\over 2}\sum_{i=1}^n\sum_{j=1}^\nu
P_{ij,q}^{(m)}(a_iy_j+x_ib_j)\right)e^{(m)}_q.$$
$$\|L_a(z)\|={1\over 2}\sum_{p=1}^n\left|\sum_{i=1}^n\sum_{j=1}^\nu
P_{ij,p}^{(f)}(a_iy_j+x_ib_j)\right|+{1\over 2}
\sum_{q=1}^\nu\left|\sum_{i=1}^n\sum_{j=1}^\nu
P_{ij,q}^{(m)}(a_iy_j+x_ib_j)\right|\leq$$
$${1\over 2}\sum_{p=1}^n\left(\sum_{i=1}^n\sum_{j=1}^\nu
P_{ij,p}^{(f)}\left|a_iy_j\right|+ \sum_{i=1}^n\sum_{j=1}^\nu
P_{ij,p}^{(f)}\left|x_ib_j\right|\right)+$$ $${1\over
2}\sum_{q=1}^\nu\left(\sum_{i=1}^n\sum_{j=1}^\nu
P_{ij,q}^{(m)}\left|a_iy_j\right|+ \sum_{i=1}^n\sum_{j=1}^\nu
P_{ij,q}^{(m)}\left|x_ib_j\right|\right)$$ Using (\ref{2}) from the
last inequality we get
$$\|L_a(z)\|\leq \|a\|\|y\|+\|b\|\|x\|\leq
\max\{\|a\|,\|b\|\}\|z\|,$$ where
$$\|a\|=\sum_{i=1}^n|a_i|, \ \ \|b\|=\sum_{j=1}^\nu |b_j|, \ \
\|x\|=\sum_{i=1}^n|x_i|, \ \ \|y\|=\sum_{j=1}^\nu |y_j|.$$ Thus
$L_a$ is bounded for any fixed $a\in\B$. \endproof

\begin{thm}\label{t3}
The algebra $\B$ is complete as a normed linear space. That is, $\B$
is a Banach space.
\end{thm}
\proof If $\{z^k\}$ converges then its all coordinates also converge
and therefore limit of $z^k$ also will be an element of $\B$. This
completes the proof.
\endproof
\begin{cor}
The algebra $\B$ is a non associative Banach algebra.
\end{cor}

\section{The derivations of $\mathcal B$}

There are many papers devoted on the subject of derivations of
genetic algebras (see e.g. \cite{c,go,ho,t}).  In \cite{ho} an
explanation of the genetic meaning of a derivation of a genetic
algebra is given. For any algebra, it is known that the space of its
derivations is a Lie algebra. The Lie algebra of derivations of a
given algebra is an important tool for studying its structure,
particularly in the non-associative case; so this is a natural
development.

 In this section we describe the set of all derivations
of $\mathcal B$. Let $D\in \Der(\mathcal B)$ and suppose

$$D(e_k^{(f)})=\sum_{i=1}^nd_{ki}^{ff}e_i^{(f)}+\sum_{l=1}^\nu
d_{kl}^{fm}e_l^{(m)}, \ \ k=1,\dots,n.$$

$$D(e_k^{(m)})=\sum_{i=1}^nd_{ki}^{mf}e_i^{(f)}+\sum_{l=1}^\nu
d_{kl}^{mm}e_l^{(m)}, \ \ k=1,\dots,\nu.$$ By the definition of
derivation $D(zt)=D(z)t+zD(t)$, we have
$$D(e_k^{(f)}e_j^{(f)})=D(e_k^{(f)})e_j^{(f)}+e_k^{(f)}D(e_j^{(f)})=$$
$$\left(\sum_{i=1}^nd_{ki}^{ff}e_i^{(f)}+\sum_{l=1}^\nu
d_{kl}^{fm}e_l^{(m)}\right)e_j^{(f)}+e_k^{(f)}\left(\sum_{i=1}^nd_{ji}^{ff}e_i^{(f)}+\sum_{l=1}^\nu
d_{jl}^{fm}e_l^{(m)}\right)=$$
$${1\over
2}\sum_{i=1}^n\left(\sum_{l=1}^\nu\left(d_{kl}^{fm}P^{(f)}_{jl,i}+d_{jl}^{fm}P^{(f)}_{kl,i}\right)\right)e_i^{(f)}+$$
$$ {1\over
2}\sum_{q=1}^\nu\left(\sum_{l=1}^\nu\left(d_{kl}^{fm}P^{(m)}_{jl,q}+d_{jl}^{fm}P^{(m)}_{kl,q}\right)\right)e_q^{(m)}=0
$$
Consequently,
\begin{equation}\label{9}\begin{array}{ll}
\sum_{l=1}^\nu\left(d_{kl}^{fm}P^{(f)}_{jl,i}+d_{jl}^{fm}P^{(f)}_{kl,i}\right)=0, \ \ i,j,k=1,\dots,n; \\[3mm]
\sum_{l=1}^\nu\left(d_{kl}^{fm}P^{(m)}_{jl,q}+d_{jl}^{fm}P^{(m)}_{kl,q}\right)=0,
\ \ j,k=1,\dots,n; \ q=1,\dots,\nu.
\end{array}
\end{equation}
Similarly, from
$D(e_k^{(m)}e_j^{(m)})=D(e_k^{(m)})e_j^{(m)}+e_k^{(m)}D(e_j^{(m)})=0$
we obtain
\begin{equation}\label{10}\begin{array}{ll}
\sum_{i=1}^n\left(d_{ki}^{mf}P^{(f)}_{ij,s}+d_{ji}^{mf}P^{(f)}_{ik,s}\right)=0, \ \ s=1,\dots,n; \ j,k=1,\dots,\nu; \\[3mm]
\sum_{i=1}^n\left(d_{ki}^{mf}P^{(m)}_{ij,q}+d_{ji}^{mf}P^{(m)}_{ik,q}\right)=0,
\ \ j,k,q=1,\dots,\nu.
\end{array}
\end{equation}
The equality
$D(e_i^{(f)}e_j^{(m)})=D(e_i^{(f)})e_j^{(m)}+e_i^{(f)}D(e_j^{(m)})$
gives the following conditions
\begin{equation}\label{11}\begin{array}{llll}
\sum_{p=1}^n\left(d_{ip}^{ff}P^{(f)}_{pj,s}-d_{ps}^{ff}P^{(f)}_{ij,p}\right)+\\[2mm]
\sum_{q=1}^\nu\left(d_{jq}^{mm}P^{(f)}_{iq,s}-d_{qs}^{mf}P^{(m)}_{ij,q}\right)=0, \ \ i,s=1,\dots,n; \ j=1,\dots,\nu; \\[4mm]
\sum_{p=1}^n\left(d_{ip}^{ff}P^{(m)}_{pj,t}-d_{pt}^{fm}P^{(f)}_{ij,p}\right)+\\[2mm]
\sum_{q=1}^\nu\left(d_{jq}^{mm}P^{(m)}_{iq,t}-d_{qt}^{mm}P^{(m)}_{ij,q}\right)=0,
\ \ i=1,\dots,n; \ j,t=1,\dots,\nu.
\end{array}
\end{equation}
Since for any $z=(x,y)\in \mathcal B$ we have $
D(z)=\sum_{i=1}^nx_iD(e^{(f)}_i)+\sum_{j=1}^\nu y_jD(e^{(m)}_j)$,
$D$ is uniquely defined by the matrix $\mathcal D={\mathcal
D}(D)=\left(d^{ff}_{ij}, d^{fm}_{ip}, d^{mf}_{pi},
d^{mm}_{pq}\right)_{{i,j=1,\dots,n\atop p,q=1,\dots,\nu}}.$

Hence,

$$\Der(\mathcal B)=\left\{D: {\mathcal D}(D) \ \ \mbox{satisfies} \ \
(\ref{9}), (\ref{10}), (\ref{11})\right\}.$$

{\bf Example.} Consider the case $n=1$, $\nu=2$. In this case
$P^{(f)}_{11,1}=P^{(f)}_{12,1}=1$. Denote $P^{(m)}_{11,1}=a$,
$P^{(m)}_{11,2}=1-a$, $P^{(m)}_{12,1}=b$, $P^{(m)}_{12,2}=1-b$.
Assume $a=1-b$, $a\ne {1\over 2}$. Then using (\ref{9}), (\ref{10}),
(\ref{11}) we get
$$ D(e^{(f)}_1)=0; \ \
D(e^{(m)}_1)=\alpha\cdot(e^{(m)}_1-e^{(m)}_2); \ \
D(e^{(m)}_2)=\alpha\cdot(-e^{(m)}_1+e^{(m)}_2),$$ where $\alpha\in
\R.$ Consequently, for arbitrary element $z=(x,y_1,y_2)\in \R^{1+2}$
we have
$$ D(z)=\alpha\cdot(y_1-y_2)(e^{(m)}_1-e^{(m)}_2).$$
Thus if $n=1$, $\nu=2$ and $a=1-b$, $a\ne {1\over 2}$ then
$$\Der(\mathcal B)=\left\{D:
D(z)=\alpha(y_1-y_2)(e^{(m)}_1-e^{(m)}_2), \ z=(x,y_1,y_2)\in
\mathcal B, \alpha\in \R\right\}.$$

\section{Dynamics of the operator (\ref{3})}

Extend the operator (\ref{3}) on $\R^{n+\nu}$, i.e. consider the
operator $V:\R^{n+\nu}\to \R^{n+\nu}$, $z=(x,y)\mapsto
z'=(x',y')=V(z)$ defined as
\begin{equation}\label{d1}
x'_j= \sum_{i,k=1}^{n,\nu}P_{ik,j}^{(f)}x_iy_k; \ \ y'_l=
\sum_{i,k=1}^{n,\nu} P_{ik,l}^{(m)}x_iy_k.
\end{equation}
Consider the following linear form $X:\R^n\to \R$ ( $Y:\R^\nu\to
\R$) defined by
\begin{equation}\label{d2}
 X(x)=\sum_{i=1}^nx_i, \ \ \left(Y(y)=\sum_{k=1}^\nu y_k\right).
\end{equation}
Denote
\begin{equation}\label{d3}
 H_i=\{z=(x,y): X(x)=Y(y)=i\}, \ \ i=0,1,
\end{equation}
the product of the $i$-hyperplanes in $\R^n$ and $\R^\nu$,
respectively.

A point $z=(x,y)\in\R^{n+\nu}$ is called a {\it fixed point} (resp.
{\it zero point}) of $V$ if $V(z)=z$ (resp. $V(z)=0$).

\begin{proposition}\label{dp1}
\begin{itemize}
  \item[(1)] If $z$ is a fixed point then $z\in H_0\cup H_1$.
  \item[(2)] If $z$ is a zero point then $z\in \{z=(x,y): X(x)Y(y)=0\}$.
\end{itemize}

\end{proposition}
\proof From (\ref{d1}) and using (\ref{2})  we get
\begin{equation}\label{d4}\begin{array}{ll}
X(x')=\sum_{j=1}^nx'_j=\sum_{i,k=1}^{n,\nu}\left(\sum_{j=1}^nP_{ik,j}^{(f)}\right)x_iy_k=
X(x)Y(y)\\[2mm]
Y(y')=\sum_{l=1}^\nu y'_l=\sum_{i,k=1}^{n,\nu}\left(\sum_{l=1}^\nu
P_{ik,l}^{(m)}\right)x_iy_k= X(x)Y(y).\end{array}
\end{equation}
\begin{itemize}
  \item[(1)] If $z$ is a fixed point then $X(x')=X(x)$ and $Y(y')=Y(y)$, this
by (\ref{d4}) gives that $X(x)=Y(y)$ and $X(x)=(X(x))^2$. Hence
$X(x)=Y(y)=0$ or $1$.
 \item[(2)] If $z$ is a zero point then $X(x')=Y(y')=0$, and (\ref{d4}) gives
$X(x)Y(y)=0$.
\end{itemize}
\endproof

Let $z^{(0)}=(x^{(0)},y^{(0)})\in\R^{n+\nu}$ be an initial point.
Its {\it trajectory} is defined by
$z^{(t)}=(x^{(t)},y^{(t)})=V(z^{(t-1)})$, $t=1,2,\dots$. Denote by
$\omega(z^{(0)})$ the set of limit points of the trajectory
$\{z^{(t)}\}_{t=0}^\infty$.  If $\omega(z^{(0)})$ consists of a
single point, then the trajectory converges, and $\omega(z^{(0)})$
is a fixed point of the operator $V$.

The following theorem gives an upper estimate of the set
$\omega(z^{(0)})$

\begin{thm}\label{td1} We have
$$\omega(z^{(0)})\subset \left\{\begin{array}{lll}
H_0 \ \ \mbox{if} \ \ \left|X(x^{(0)})Y(y^{(0)})\right|<1\\[2mm]
H_1 \ \ \mbox{if} \ \ \left|X(x^{(0)})Y(y^{(0)})\right|=1\\[2mm]
H_\infty \ \ \mbox{if} \ \ \left|X(x^{(0)})Y(y^{(0)})\right|>1,
\end{array}\right.
$$
where $H_\infty=\{z=(x,y): X(x)=Y(y)=+\infty\}$.
\end{thm}
\proof Using (\ref{d4}) we get
$X(x^{(t)})=Y(y^{(t)})=X(x^{(t-1)})Y(y^{(t-1)})$ for any
$t=1,2,\dots$. Iterating this recurrent equation we obtain
\begin{equation}\label{d5}
X(x^{(t)})=Y(y^{(t)})=\left(X(x^{(0)})Y(y^{(0)})\right)^{2^{t-1}}.
\end{equation}
Consequently,
$$\lim_{t\to\infty}X(x^{(t)})=\lim_{t\to\infty}Y(y^{(t)})=\left\{\begin{array}{lll}
0 \ \ \mbox{if} \ \ \left| X(x^{(0)})Y(y^{(0)})\right|<1\\[2mm]
1 \ \ \mbox{if} \ \  \left|X(x^{(0)})Y(y^{(0)})\right|=1\\[2mm]
+\infty \ \ \mbox{if} \ \ \left|X(x^{(0)})Y(y^{(0)})\right|>1.
\end{array}\right.
$$
\endproof

\section{A special case of an EABP}

In this section we consider a special case of an EABP giving an
additional condition on heredity coefficients (\ref{2}), i.e.
consider the coefficients as follows
\begin{equation}\label{v1}P^{(f)}_{ik,j}=\left\{\begin{array}{lll}
a_{ij} \ \ \mbox{if} \ \ k=1 \\[2mm]
1 \ \ \ \ \mbox{if} \ \ k\ne 1, j=1 \\[2mm]
0 \ \ \ \ \mbox{if} \ \ k\ne 1, j\ne 1.
\end{array}\right. ,
\ \  P^{(m)}_{ik,l}=\left\{
\begin{array}{lll}
b_{kl} \ \ \mbox{if} \ \ i=1 \\[2mm]
1 \ \ \ \ \mbox{if} \ \ i\ne 1, l=1 \\[2mm]
0 \ \ \ \ \mbox{if} \ \ i\ne 1, l\ne 1.
\end{array} \right.
\end{equation}
The matrices $A=(a_{ij})$ and $B=(b_{kl})$ by (\ref{2}) satisfy the
following conditions
\begin{equation}\label{v}
 a_{ij}\geq 0, \ \ \sum_{j=1}^na_{ij}=1, \ \
i=1,\dots,n; \ \ b_{kl}\geq 0, \ \ \sum_{l=1}^\nu b_{kl}=1, \ \
k=1,\dots,\nu,\end{equation} i.e. both matrices are stochastic.

\begin{rk} The condition (\ref{v1}) \ is taken to simplify our computations;
in this way we reduced both cubic matrices
$\left(P^{(f)}_{ik,j}\right)$, $ \left(P^{(m)}_{ik,l}\right)$ to
quadratic matrices. But that condition has very clear biological
treatment: the type ``1'' of females and the type ``1'' of males
have preference, i.e. any type of female (male) can be born if its
father (mother) has type ``1''. If the father (mother) has type $\ne
1$ then only type ``1'' female (male) can be born.
\end{rk}
Under condition (\ref{v1}) the multiplication (\ref{4}) became as

\begin{equation}\label{v2}
\begin{array}{ll}
e^{(f)}_ie^{(m)}_k = e^{(m)}_ke^{(f)}_i={1\over
2}\left\{\begin{array}{llll} \sum_{j=1}^na_{1j}e^{(f)}_j+
\sum_{l=1}^{\nu}b_{1l}e^{(m)}_l, \ \ \mbox{if} \ \ i=1,
k=1\\[2mm]
\sum_{j=1}^na_{ij}e^{(f)}_j+ e^{(m)}_1, \ \ \mbox{if} \ \
i\ne 1, k=1\\[2mm]
e^{(f)}_1+\sum_{l=1}^\nu b_{kl}e^{(m)}_l, \ \ \mbox{if} \ \
i =1, k\ne 1\\[2mm]
e^{(f)}_1+e^{(m)}_1, \ \ \mbox{if} \ \ i\ne 1, k\ne 1.
\end{array}\right.\\[5mm]
e^{(f)}_ie^{(f)}_j =0,\ \ i,j=1,\dots,n; \ \ e^{(m)}_ke^{(m)}_l =0,
\ \ k,l=1,\dots,\nu. \end{array} \end{equation}

Operator (\ref{3}) has the following form

\begin{equation}\label{v3}\begin{array}{llll}
x'_1=\sum_{i=1}^n\left(a_{i1}y_1+\sum_{k=2}^\nu y_k\right)x_i\\[2mm]
x'_j=y_1\sum_{i=1}^na_{ij}x_i, \ \ j\ne 1\\[4mm]
y'_1=\sum_{k=1}^\nu\left(b_{k1}x_1+\sum_{i=2}^n x_i\right)y_k\\[2mm]
y'_l=x_1\sum_{k=1}^\nu b_{kl}y_k, \ \ l\ne 1
\end{array}
\end{equation}
Denote by $\B_1$ the EABP defined by the multiplication table
(\ref{v2}).

\subsection{Idempotent elements of $\B_1$}

A element $z\in \B$ is called {\it idempotent} if $z^2=z$; such
points of an EABP are especially important, because they are the
fixed points of the evolution map $V$, i.e. $V(z)=z$. We denote by
${\mathcal Id}(\B)$ the idempotent elements of an algebra $\B$.
Clearly, $0\in {\mathcal Id}(\B)$ and this set is an algebraic
variety. By Proposition \ref{dp1} we have ${\mathcal Id}(\B)\subset
H_0\cup H_1$. In this subsection we shall describe idempotent
elements of $\B_1$.  First we describe the set of idempotent
elements which belong in $H_0\cap {\mathcal Id}(\B_1)$ and after
that we shall describe elements of $H_1\cap {\mathcal Id}(\B_1)$.

{\bf Idempotents in $H_0$.}

 Using (\ref{v3}) and the condition that $z\in H_0$, from the equation
$V(z)=z^2=z$ we obtain
\begin{equation}\label{v4}\begin{array}{ll}
x_j=y_1\sum_{i=2}^n(a_{ij}-a_{1j})x_i, \ \ j=2,\dots,n;\\[2mm]
y_l=x_1\sum_{k=2}^\nu (b_{kl}-b_{1l})y_k, \ \ l=2,\dots,\nu.
\end{array}
\end{equation}

{\it Case} $x_1y_1=0$: If $x_1=0$ then $y_2=\dots=y_\nu=0$,
consequently, $y_1=0$. Similarly, if $y_1=0$ we get
$x_1=\dots=x_n=0$. Hence in case $x_1y_1=0$ we have unique
idempotent $z=0$.

{\it Case} $x_1y_1\ne 0$: Consider matrices
$C_y=(c_{ij})_{i,j=2,\dots,n}$ and $D_x=(d_{kl})_{k,l=2,\dots,\nu}$
such that
$$c_{ij}=\left\{\begin{array}{ll}
(a_{ij}-a_{1j})y \ \ \ \ \mbox{if} \ \ i\ne j \\[2mm]
(a_{ij}-a_{1j})y-1 \ \ \mbox{if} \ \ i=j
\end{array}\right. ,
\ \  d_{kl}=\left\{
\begin{array}{ll}
(b_{kl}-b_{1l})x \ \ \ \ \mbox{if} \ \ k\ne l \\[2mm]
(b_{kl}-b_{1l})x-1 \ \ \mbox{if} \ \ k=l
\end{array} \right.
$$ Then equation (\ref{v4}) can be written as
\begin{equation}\label{v5} C_{y_1}x=0, \ \ D_{x_1}y=0,\end{equation}
where $x=(x_2,\dots,x_n)$,  $y=(y_2,\dots,y_{\nu})$.

Consider now $x_1\in \R\setminus\{0\}$ as a parameter, then equation
$D_{x_1}y=0$ has a unique solution $y=0$ if $\det(D_{x_1})\ne 0$,
which gives $y_1=0$, i.e. this is a contradiction to the assumption
that $y_1\ne 0$.

If $\det(D_{x_1})= 0$ then we fix a solution $x_1=x_1^*\ne 0$ of the
equation $\det(D_{x_1})= 0$. In this case there are infinitely many
solutions $y^*=(y^*_2,\dots,y^*_\nu)$ of $D_{x_1}y= 0$. Substituting
the solution $y^*_1=-\sum_ {k=2}^\nu y^*_k$ in $C_{y_1}x=0$, we get
\begin{equation}\label{v6}
x_j=y^*_1\sum_{i=2}^n(a_{ij}-a_{1j})x_i=0,
j=2,\dots,n.\end{equation} This system has a unique solution
$x_2=\dots=x_n=0$ if $\det(C_{y^*_1})\ne 0$ but in this case we get
$x_1=0$ which is a contradiction to the assumption that $x_1\ne 0$.
If $\det(C_{y^*_1})=0$ then we have infinitely many solutions
$x^*=(x^*_2,...,x^*_n)$.

Hence we have proved the following
\begin{proposition}\label{p1} We have
$$H_0\cap {\mathcal Id}(\B_1)=\{0\}\cup$$ $$
\{((x_1^*,...,x_n^*),(y_1^*,...,y_\nu^*)): C_{y^*_1}x^*=0, \ \
D_{x^*_1}y^*=0, \det(D_{x^*_1})=\det(C_{y^*_1})=0\}. $$
\end{proposition}

{\bf Idempotents in $H_1$.} Using (\ref{v3}) and the condition that
$z\in H_1$, the equation $V(z)=z^2=z$ can be written as

\begin{equation}\label{v3x}\begin{array}{ll}
x_1=\sum_{i=1}^n\left((a_{i1}-1)y_1+1\right)x_i;\\[2mm]
x_j=y_1\sum_{i=1}^na_{ij}x_i, \ \ j\ne 1.
\end{array}
\end{equation}
\begin{equation}\label{v3y}\begin{array}{ll}
y_1=\sum_{k=1}^\nu\left((b_{k1}-1)x_1+1\right)y_k;\\[2mm]
y_l=x_1\sum_{k=1}^\nu b_{kl}y_k, \ \ l\ne 1.
\end{array}
\end{equation}
Consider several cases.

{\bf Case} $x_1=y_1=0$.  In this case we get $z=0$ which is not in
$H_1$.

{\bf Case} $x_1\ne 0,$ $y_1=0$. In this case $x_1=1$, $x_j=0, \
j=2,\dots,n$. Substituting $x_1=1$ in the system of equations
(\ref{v3y}) we get
\begin{equation}\label{v3z}\begin{array}{ll}
0=\sum_{k=2}^\nu b_{k1}y_k\\[2mm]
y_l=\sum_{k=2}^\nu b_{kl}y_k, \ \ l=2,\dots,\nu.
\end{array}
\end{equation}

Denote $B_2=\left(b'_{kl}\right)_{k,l=2,...,\nu}$, with
$$ b'_{kl}=\left\{
\begin{array}{ll}
b_{kl}-1 \ \ \ \ \mbox{if} \ \ k=l \\[2mm]
b_{kl} \ \ \mbox{if} \ \ k\ne l.
\end{array} \right.
$$

If $\det(B_2)\ne 0$ then $B_2y=0$ gives $y_1=\dots=y_\nu=0$ but this
does not satisfy $\sum_{k=2}^\nu y_k=1$. If $\det(B_2)= 0$ then
$B_2y=0$ has infinitely many solutions $y=(y_2,\dots,y_\nu)$. We
must take these solutions which satisfy $\sum_{k=2}^\nu
b_{k1}y_k=0$. So in this case we have the following idempotent
elements (which belong to $H_1$)
$$I_0=\left\{\begin{array}{ll}
\{((1,0,\dots,0),(0,y_2,\dots,y_\nu)): B_2y=0,\sum_{k=2}^\nu
b_{k1}y_k=0\} \  \mbox{if} \ \det(B_2)= 0\\[2mm]
\quad \emptyset \hspace{9cm} \mbox{if}
 \ \det(B_2)\ne 0.
\end{array}\right.$$

{\bf Case} $x_1=0$, $y_1\ne 0$. This case is similar to the previous
case. Denote $A_2=\left(a'_{ij}\right)_{i,j=2,...,n}$, with
$$ a'_{ij}=\left\{
\begin{array}{ll}
a_{ij}-1 \ \ \ \ \mbox{if} \ \ i=j \\[2mm]
a_{ij} \ \ \mbox{if} \ \ i\ne j.
\end{array} \right.
$$

Then we have the following idempotent elements (which belong to
$H_1$):
$$I_1=\left\{\begin{array}{ll}
\{((0,x_1,\dots,x_n),(1,0,\dots,0)): A_2x=0,\sum_{i=2}^n
a_{i1}x_i=0\} \  \mbox{if} \ \det(A_2)= 0\\[2mm]
\quad \emptyset \hspace{9cm} \mbox{if} \ \det(A_2)\ne 0.
\end{array}\right.$$

{\bf Case} $x_1\ne 0$, $y_1\ne 0$. Here, to avoid many special cases
we assume that $a_{11}\ne 1$.  Take $y_1$ as a parameter and solve
the system (\ref{v3x}) which is equivalent to the following
\begin{equation}\label{v14d}
\sum_{i=2}^n\left({a_{1j}((a_{i1}-1)y_1+1)\over
1-a_{11}}+a_{ij}y_1\right)x_i=x_j, \ \ j=2,\dots,n.
\end{equation}
Denote $U_y=(u_{ij})_{i,j=2,\dots,n}$ with
$$ u_{ij}=\left\{
\begin{array}{ll}
{a_{1j}((a_{i1}-1)y+1)\over
1-a_{11}}+a_{ij}y  \ \ \mbox{if} \ \ i\ne j \\[2mm]
{a_{1j}((a_{i1}-1)y+1)\over 1-a_{11}}+a_{ij}y-1 \ \ \mbox{if} \ \ i=
j
\end{array} \right.
$$

{\bf Subcase} $\det(U_y)\ne 0$. In this case we have
$x_2=\dots=x_n=0$, consequently, $x_1=1$. For $x_1=1$ from the
equation (\ref{v3y}) we get $B_1y=0$ where
$B_1=(b^*_{kl})_{k,l=1,...,\nu}$, with
$$ b^*_{kl}=\left\{
\begin{array}{ll}
b_{kl}-1 \ \ \ \ \mbox{if} \ \ k=l \\[2mm]
b_{kl} \ \ \mbox{if} \ \ k\ne l.
\end{array} \right.
$$
So in this case we have the following set of idempotent elements of
$H_1$.
$$I_2=\left\{\begin{array}{ll}
\{((1,0,\dots,0),(y_1,\dots,y_\nu)): B_1y=0, \det(U_y)\ne 0\} \  \mbox{if} \ \det(B_1)= 0\\[2mm]
\quad \emptyset \hspace{8cm} \mbox{if} \ \det(B_1)\ne 0.
\end{array}\right.$$

{\bf Subcase} $\det(U_y)= 0$. In this case we fix a solution
$y_1=y_1^*$ of $\det(U_y)= 0$. We have infinitely many solutions
$(x^*_2,...,x^*_n)$. Denote $C_x=(c^*_{kl})_{k,l=1,...,\nu}$, with
$$ c^*_{kl}=\left\{
\begin{array}{llll}
(b_{11}-1)x \ \ \ \ \mbox{if} \ \ k=l=1; \\[2mm]
(b_{k1}-1)x+1 \ \ \mbox{if} \ \ k\ne 1, l=1;\\[2mm]
b_{ll}x-1 \ \ \ \ \mbox{if} \ \ k=l\ne 1; \\[2mm]
b_{k1}x \ \ \mbox{if} \ \ k\ne l, l=2,\dots,\nu.\\[2mm]
\end{array} \right.
$$
Now one has the following set of idempotent elements
$$I_3=\{((x^*_1,\dots,x^*_n),(y^*_1,\dots,y^*_\nu)): C_{x_1^*}y^*=0,
U_{y_1^*}x^*=0,$$ $$ \det(U_{y_1^*})=\det(C_{x_1^*})=0\}.$$

Thus we have the following
\begin{thm}\label{t4d} If $a_{11}\ne 1$ then the full set of the idempotent elements which belong to $H_1$ is
$${\mathcal Id}(\B_1)\cap H_1 =I_0\cup I_1\cup I_2\cup I_3.
 $$
\end{thm}

  \begin{rk}
By Proposition \ref{dp1} we can conclude that Proposition \ref{p1}
and Theorem \ref{t4d} give the full set ${\mathcal Id}(\B_1)$. But
conditions described in the Proposition and Theorem are complicated
in general. In the sequel of this subsection we shall consider
additional conditions on $(a_{ij})$, $(b_{ij})$ and under these
conditions we explicitly describe the set ${\mathcal Id}(\B_1)$.
\end{rk}

Now we consider a particular case and describe the full set of
idempotent elements of $\B_1$.

If we consider the case
\begin{equation}\label{v7} a_{ij}=\left\{\begin{array}{ll}
1, \ \ i=j\\[2mm]
0, \ \ i\ne j,
\end{array}\right.\ \
 b_{kl}=\left\{\begin{array}{ll}
1, \ \ k=l\\[2mm]
0, \ \ k\ne l
\end{array}\right.
\end{equation}
then the following is true
\begin{proposition}\label{p3} If the condition (\ref{v7}) is
satisfied then
$${\mathcal Id}(\B_1)=\{0\}\cup\left\{((1,x_2,\dots,x_n),(1,y_2,\dots,y_\nu)):\sum_{i=2}^n x_i=\sum_{k=2}^\nu y_k=0\right\}\cup$$ $$
\left\{((1,x_2,\dots,x_n),(1,y_2,\dots,y_\nu)):\sum_{i=2}^n
x_i=\sum_{k=2}^\nu y_k=-1 \right\}\cup$$
$$\left\{((1,0,\dots,0),(y_1,\dots,y_\nu)):\sum_{k=1}^\nu y_k=1, y_1\ne 1\right\}\cup$$
$$\left\{((x_1,\dots,x_n),(1,0,\dots,0)): \sum_{i=1}^n x_i=1, x_1\ne 1\right\}.$$
\end{proposition}
\proof The equation $z^2=z$, for $z=(x,y)\in \R^{n+\nu}$, using
condition (\ref{v7}) can be written as
\begin{equation}\label{v10}\begin{array}{llll}
(1-y_1)x_1=\sum_{k=2}^\nu y_k\sum_{i=1}^nx_i,\\[2mm]
(1-y_1)x_j=0, \ \ j=2,\dots,n\\[3mm]
(1-x_1)y_1=\sum_{i=2}^nx_i\sum_{k=1}^\nu y_k,\\[2mm]
(1-x_1)y_l=0, \ \ l=2,\dots,\nu
\end{array}
\end{equation}
The simple analysis of the system (\ref{v10}) gives the set of all
idempotents.
\endproof
\subsection{Absolute nilpotent elements of $\B_1$}
The element $z$ is called an {\it absolute nilpotent} if $z^2=0$,
i.e. $z$ is zero-point of the evolution operator $V$. For $\B_1$ the
equation $z^2=0$ is equivalent to the following system of quadratic
equations
\begin{equation}\label{n1}\begin{array}{llll}
\sum_{i=1}^n\left(a_{i1}y_1+\sum_{k=2}^{\nu}y_k\right)x_i=0;\\[2mm]
y_1\sum_{i=1}^na_{ij}x_i=0, \ \ j=2,\dots,n;\\[4mm]
\sum_{k=1}^\nu\left(b_{k1}x_1+\sum_{i=2}^nx_i\right)y_k=0;\\[2mm]
x_1\sum_{k=1}^\nu b_{kl}y_k=0, \ \ l=2,\dots,\nu.
\end{array}
\end{equation}

By Proposition \ref{dp1} we have
$\sum_{i=1}^nx_i\sum_{k=1}^{\nu}y_k=0$. There are the following
three cases.
\begin{itemize}
\item {\it Case} $\sum_{i=1}^nx_i=\sum_{k=1}^{\nu}y_k=0$.
In this case, from  the system of quadratic equations (\ref{n1}) we
get
\begin{equation}\label{n2}\begin{array}{ll}
y_1\sum_{i=1}^na_{ij}x_i=0, \ \ j=1,\dots,n;\\[3mm]
x_1\sum_{k=1}^\nu b_{kl}y_k=0, \ \ l=1,\dots,\nu.
\end{array}
\end{equation}

{\it Subcase} $y_1=x_1=0$. In this case one easily gets the
following set of solutions
$${\mathcal N}^0_{00}=\left\{((0,x_2,\dots,x_n),(0,y_2,\dots,y_\nu)):
\sum_{i=2}^nx_i=0, \sum_{k=2}^\nu y_k=0\right\}.$$

{\it Subcase} $x_1=0, y_1\ne 0$. Denote
$A_0=(a_{ij})_{i,j=2,\dots,n}$. It is easy to see that the set of
solutions is as follows

$${\mathcal N}^0_{01}= \left\{\begin{array}{ll}
{\mathcal N}^1 \ \ \mbox{if} \ \
{\rm det}(A_0)\ne 0;\\[2mm]
{\mathcal N}^0 \ \ \mbox{if} \ \ {\rm det}(A_0)=
0,\end{array}\right.$$ where
$${\mathcal N}^1=\left\{((0,\dots,0),(y_1,\dots,y_\nu)): y_1\ne 0,
\sum_{k=1}^\nu y_k=0\right\},$$ $${\mathcal
N}^0=\left\{((0,x_2,\dots,x_n),(y_1,\dots,y_{\nu})):
\sum_{i=2}^nx_i=\sum_{i=2}^na_{i1}x_i=0, \right.$$
$$\left.A_0x=0, y_1\ne
0, \sum_{k=1}^\nu y_k=0\right\}.
$$

{\it Subcase} $x_1\ne 0, y_1=0$. Denote
$B_0=(b_{kl})_{k,l=2,\dots,\nu}$. In this case the set of solutions
is
$${\mathcal N}^0_{10}= \left\{\begin{array}{ll}
{\mathcal N}^2 \ \ \mbox{if} \ \
{\rm det}(B_0)\ne 0;\\[2mm]
{\mathcal N}^3 \ \ \mbox{if} \ \ {\rm det}(B_0)=
0,\end{array}\right.$$ where
$${\mathcal N}^2=\left\{((x_1,\dots,x_n),(0,\dots,0)): x_1\ne 0,
\sum_{i=1}^n x_i=0\right\},$$ $${\mathcal
N}^3=\left\{((x_1,\dots,x_n),(0,y_2,\dots,y_{\nu})): x_1\ne 0,
\sum_{i=1}^n x_i=0,\right.$$ $$ \left. B_0y=0,  \sum_{k=2}^\nu y_k=
\sum_{k=2}^\nu b_{k1}y_k =0\right\}.
$$

{\it Subcase} $x_1\ne 0, y_1\ne 0$. In this case it is easy to get
the following

$${\mathcal N}^0_{11}= \left\{\begin{array}{ll}
\emptyset \ \ \mbox{if} \ \
{\rm det}(A)\ne 0 \ \ \mbox{or} \ \ {\rm det}(B)\ne 0;\\[2mm]
{\mathcal N}^4 \ \ \mbox{if} \ \ {\rm det}(A)={\rm
det}(B)=0,\end{array}\right.$$ where
$${\mathcal N}^4=\left\{((x_1,\dots,x_n),(y_1,\dots,y_\nu)): Ax=0, By=0\right\}.$$

\item {\it Case} $\sum_{i=1}^nx_i=0, \sum_{k=1}^{\nu}y_k\ne 0$.
In this case the system of quadratic equations (\ref{n1}) can be
written as
\begin{equation}\label{n3}\begin{array}{llll}
\sum_{i=1}^n\left(a_{i1}y_1+\sum_{k=2}^{\nu}y_k\right)x_i=0\\[2mm]
y_1\sum_{i=1}^na_{ij}x_i=0, \ \ j=2,\dots,n\\[4mm]
x_1\sum_{k=1}^\nu\left(b_{k1}-1\right)y_k=0\\[2mm]
x_1\sum_{k=1}^\nu b_{kl}y_k=0, \ \ l=2,\dots,\nu.
\end{array}
\end{equation}

{\it Subcase} $y_1=x_1=0$. In this case we have the following set of
absolute nilpotents

$${\mathcal N}^1_{00}=\left\{((0,x_2,\dots,x_n),(0,y_2,\dots,y_\nu)): \sum_{i=2}^nx_i=0, \sum_{k=2}^\nu y_k\ne 0\right\}.$$

{\it Subcase} $x_1=0, y_1\ne 0$. From (\ref{n3}) we get
\begin{equation}\label{n4}\begin{array}{ll}
\sum_{i=2}^n\left(a_{i1}y_1+\sum_{k=2}^{\nu}y_k\right)x_i=0\\[2mm]
\sum_{i=1}^na_{ij}x_i=0, \ \ j=2,\dots,n.\\[4mm]
\end{array}
\end{equation}
This has the following set of solutions

$${\mathcal N}^1_{01}= \left\{\begin{array}{ll}
{\mathcal N}^5 \ \ \mbox{if} \ \
{\rm det}(A_0)\ne 0;\\[2mm]
{\mathcal N}^6 \ \ \mbox{if} \ \ {\rm det}(A_0)=
0,\end{array}\right.$$ where
$${\mathcal N}^5=\left\{((0,\dots,0),(y_1,\dots,y_\nu)): y_1\ne 0,
\sum_{k=1}^\nu y_k\ne 0\right\},$$ $${\mathcal
N}^6=\left\{((0,x_2,\dots,x_n),(y_1,\dots,y_{\nu})):
\sum_{i=2}^nx_i=\sum_{i=2}^n\left(a_{i1}y_1+\sum_{k=2}^\nu
y_k\right)x_i=0,\right.$$
$$ \left.A_0x=0,
y_1\ne 0, \sum_{k=1}^\nu y_k\ne 0\right\}.
$$

{\it Subcase} $x_1\ne 0, y_1=0$. In this case from (\ref{n3}) we get
\begin{equation}\label{n4}\begin{array}{llll}
\sum_{k=2}^\nu\left(b_{k1}-1\right)y_k=0\\[2mm]
\sum_{k=2}^\nu b_{kl}y_k=0, \ \ l=2,\dots,\nu.
\end{array}
\end{equation}

If $\det(B_0)\ne 0$ we get $y_2=\dots=y_\nu=0$, this is impossible
since we have condition $\sum_{k=1}^{\nu}y_k\ne 0$. For $\det(B_0)=
0$ the set of absolute nilpotents will be

$${\mathcal N}^1_{10}=\left\{((x_1,\dots,x_n),(0,y_2,\dots,y_\nu)): x_1\ne 0,
\sum_{i=1}^n x_i=0, \right.$$ $$ \left.B_0y=0,
\sum_{k=2}^{\nu}y_k\ne 0,
\sum_{k=2}^\nu\left(b_{k1}-1\right)y_k=0\right\}.$$

{\it Subcase} $x_1\ne 0, y_1\ne 0$. Denote ${\bf A}=({\bf
a}_{ij})_{i,j=1,\dots,n}$, with
$${\bf a}_{ij}=\left\{
\begin{array}{ll}
1 \ \ \ \ \mbox{if} \ \ j=1; \\[2mm]
a_{ij} \ \ \mbox{if} \ \ j\ne 1.
\end{array} \right.
$$
and  ${\bf B}=({\bf b}_{kl})_{k,l=1,\dots,\nu}$, with
$${\bf b}_{kl}=\left\{
\begin{array}{ll}
b_{k1}-1 \ \ \ \ \mbox{if} \ \ l=1 \\[2mm]
b_{kl} \ \ \mbox{if} \ \ l\ne 1.
\end{array} \right.
$$
In this case we have
$${\mathcal N}^1_{11}= \left\{\begin{array}{ll}
\emptyset \ \ \mbox{if} \ \
{\rm det}({\bf A})\ne 0 \ \ \mbox{or} \ \ {\rm det}({\bf B})\ne 0;\\[2mm]
{\mathcal N}^9 \ \ \mbox{if} \ \ {\rm det}({\bf A})={\rm det}({\bf
B})=0,\end{array}\right.$$ where
$${\mathcal N}^9=\left\{((x_1,\dots,x_n),(y_1,\dots,y_\nu)): {\bf A}x=0, {\bf B}y=0,\right.$$
$$\left. \sum_{i=1}^n\left(a_{i1}y_1+\sum_{k=2}^\nu
y_k\right)x_i=0, \sum_{k=1}^\nu y_k\ne 0 \right\}.$$

\item {\it Case} $\sum_{i=1}^nx_i\ne 0, \sum_{k=1}^{\nu}y_k=0$.
This case is similar to the previous case, to get the set of
nilpotents, one has to change sets ${\mathcal N}^1_{ij}$, $i,j=0,1$
as follows. Replace $x$ and $y$, also rearrange parameters $a_{ij}$
with $b_{kl}$. Denote resulting sets by ${\mathcal N}^2_{ij}$,
$i,j=0,1$ respectively.
\end{itemize}

Thus we have proved the following theorem

\begin{thm}\label{t4} The full set ${\mathcal N}$ of absolute nilpotent
elements of the algebra $\B_1$ is
$${\mathcal N}=\bigcup_{i=0}^2\left({\mathcal N}_{00}^i\cup{\mathcal N}_{01}^i\cup{\mathcal N}_{10}^i\cup{\mathcal N}^i_{11}\right).
 $$
\end{thm}

\section*{ Acknowledgements}

 The first author was supported by Ministerio
de Ciencia e Innovaci\'on (European FEDER support included), grant
MTM2009-14464-C02-01, and by Xunta de Galicia, grant Incite09 207
215 PR. The second author thanks the Department of Algebra, University of
Santiago de Compostela, Spain,  for providing financial support of
his visit to the Department (February-April 2010).

{}

\begin{thebibliography}{00}

\bibitem{ber} S. N. Bernstein, \emph{The solution of a mathematical problem related to the
theory of heredity}, Uchen. Zapiski Nauchno-Issled. Kafedry Ukr.
Otd. Matem. 1, 83--115 (1924).

\bibitem{c} R. Costa, \emph{On the derivations of gametic algebras for polyploidy
with multiple alleles},  Bol. Soc. Brasil. Mat.  13(2), 69--81
(1982).

\bibitem{e1} I.M.H. Etherington, \emph{Genetic algebras}, Proc. Roy. Soc. Edinburgh. 59, 242--258 (1939).

\bibitem{e2} I.M.H. Etherington, \emph{Duplication of linear algebras}, Proc. Edinburgh Math. Soc.2, 6, 222--230 (1941).

\bibitem{e3} I.M.H. Etherington, \emph{Non-associative algebra and the symbolism of genetics},
Proc. Roy. Soc. Edinburgh. 61, 24--42 (1941).

\bibitem{gr} N. N. Ganikhodjaev and U. A. Rozikov, \emph{On quadratic stochastic operators generated by
Gibbs distributions}, Regul. Chaotic Dyn. 11(4), 467--473 (2006).

\bibitem{gan1} R. N. Ganikhodzhaev, \emph{Quadratic stochastic operators, Lyapunov functions,
and tournaments}, Russian Acad. Sci. Sb. Math. 76(2), 489--506
(1993).

\bibitem{ge} R. N. Ganikhodzhaev and D. B. Eshmamatova, \emph{Quadratic automorphisms of a simplex and the
asymptotic behavior of their trajectories}, Vladikavkaz. Math. Zh.
8(2), 12--28 (2006).

\bibitem{gar} R. N. Ganikhodzhaev and U.A. Rozikov, \emph{Quadratic stochastic operators:
results and open problems}, arXiv:0902.4207v2 [math.DS]

\bibitem{go} H. Gonshor, \emph{Derivations in genetic algebras},  Comm. Algebra
16(8), 1525--1542 (1988).

\bibitem{h} P. Holgate,  \emph{Genetic algebras associated with sex linkage},  Proc.
Edinburgh Math. Soc.2, 17 113--120 (1970/71).

\bibitem{ho} P. Holgate, \emph{The interpretation of derivations in genetic
algebras}, Linear Algebra Appl.  85, 75--79  (1987).

\bibitem{ly} Y.I. Lyubich, \emph{Mathematical structures in population
genetics}, Springer-Verlag, Berlin, 1992.

\bibitem{m} M.L. Reed, \emph{Algebraic structure of genetic inheritance},  Bull. Amer. Math.
Soc. (N.S.)  34(2), 107--130  (1997).

\bibitem{rz}  U.A. Rozikov and U.U. Zhamilov, \emph{On $F$-quadratic stochastic operators},
 Math. Notes. 83(4), 554--559 (2008).

\bibitem{rz1}  U.A. Rozikov and U.U. Zhamilov, \emph{The dynamics of strictly non-Volterra
quadratic stochastic operators on the 2-simplex}, Sbornik: Math.
200(9), 1339--1351 (2009).

\bibitem{roz}  U.A. Rozikov and  A. Zada, \emph{On $\ell$-Volterra quadratic stochastic
operators}, Doklady Math. 79(1), 32--34 (2009).

\bibitem{rs} U. A. Rozikov and N. B. Shamsiddinov, \emph{On non-Volterra quadratic stochastic
operators generated by a product measure}, Stoch, Anal. Appl. 27(2)
353--362 (2009).

\bibitem{t} J. P. Tian, \emph{Evolution algebras and their applications},
Lecture Notes in Mathematics, 1921, Springer-Verlag, Berlin, 2008.

\bibitem{w} A. W$\ddot{o}$rz-Busekros, \emph{Algebras in genetics}, Lecture Notes in
Biomathematics, 36. Springer-Verlag, Berlin-New York, 1980.
\end{thebibliography}
\end{document}